\documentclass{amsart} 
\usepackage[utf8]{inputenc} 

\usepackage{amsmath, amsthm, amsfonts, amssymb}
\usepackage{tikz-cd}
\usepackage{hyperref}

\newtheorem{thm}{Theorem}[section]

\newtheorem{obsv}[thm]{Observation}
\theoremstyle{definition}
\newtheorem{defn}[thm]{Definition}
\newtheorem{defn/lem}[thm]{Definition/Lemma}
\newtheorem{rem}[thm]{Remark}

\newtheorem{problem}{Problem}[]

\newtheorem{que}[problem]{Question}

\newcommand{\C}{{\mathbb C}}

\newcommand{\R}{{\mathbb R}}

\newcommand{\Z}{{\mathbb Z}}

\newcommand{\HH}{{\mathcal H}}

\newcommand{\OO}{{\mathcal O}}

\newcommand{\UU}{{\mathcal U}}

\date{\today}

\title{Representing some II$_1$ factors in $L^2(\Lambda \backslash G)$}
\author{Lauren C. Ruth}

\address{Department of Mathematics, Vanderbilt University, 1326 Stevenson Center, Nashville, TN 37240, USA}
\email{lauren.c.ruth@vanderbilt.edu}

\begin{document}

\maketitle 

\begin{abstract} 
Let $G$ be $PGL(n,F)$, $n \geq 3$, $F$ a certain non-archimedean local field; or let $G$ be $PSL(2,\R) \times \cdots \times PSL(2,\R)$.  Let $\Gamma$ be a lattice in $G$, and let $( \Lambda_n )$ be a sequence of lattices in $G$ satisfying the pointwise limit multiplicity property.  In this note, we explain how the pointwise limit multiplicity property can be combined with a generalization of a theorem in \cite{ghj} to give representations of the II$_1$ factor $R \Gamma$ on a subspace of $L^2(\Lambda_i \backslash G)$ for some $\Lambda_i$ in $( \Lambda_n )$.  This extends a result in the author's dissertation \cite{ruthphd}. 
\end{abstract}

\section{Introduction} 

The main purpose of this note is to explain:  

\begin{obsv} \label{obsv} Let $G$ and $(\Lambda_n)$ be as in one of the following two cases.
\begin{itemize}
    \item[\textbf{Case 1.}] $G$ is $PGL(n,F)$, $n\geq 3$, $F$ a non-archimedean local field of characteristic $0$ and residue field of order relatively prime to $n$, and $( \Lambda_n )$ is any sequence of pairwise non-conjugate lattices in $G$; or,   
    \item[\textbf{Case 2.}] $G$ is $PSL(2,\R)^r$, $r \geq 2$, and $( \Lambda_n )$ is a sequence of lattices in $G$ with each $\Lambda_i$ of the form $PSL(2,\OO_F)$, where $\OO_F$ is the ring of integers of a real number field $F$ with $r$ real embeddings.  
\end{itemize}
Let $\Gamma$ be a lattice in $G$.  Then there exists a representation of the II$_1$ factor $R\Gamma$ on a subspace of $L^2(\Lambda_i \backslash G)$ for some $\Lambda_i$ in $( \Lambda_n )$.
\end{obsv}

The content of the proof of this observation lies in work on the limit multiplicity property carried out in \cite{7s}, \cite{gellev}, and \cite{matz}, as well as work on representations of II$_1$ factors in \cite{ghj} rooted in Atiyah's $L^2$-index theory in \cite{ati}, Atiyah and Schmid's work on discrete series representations in \cite{atischmid}, and Murray and von Neumann's foundational work \cite{roo1}.  We give the proof in Section \ref{proof}.

\begin{rem}
We could have featured groups besides $PGL(n,F)$, $n \geq 3$, and $PSL(2,\R)\times \ldots \times PSL(2,\R)$, since there are more general results on the pointwise limit multiplicity property than those we mention here.  But these two groups, along with the group $PSL(2,\R)$ in Theorem \ref{vnreppsl2r} below, already demonstrate that representations of the II$_1$ factor $R\Gamma$ on subspaces of $L^2(\Lambda \backslash G)$ can be obtained regardless of whether the ambient group $G$ is connected or totally disconnected; higher-rank or rank-1; and Property (T) or a-T-menable.  
\end{rem}

\section{Background}

First we give background on representations of $R\Gamma$, then we give background on the pointwise limit multiplicity property.  

In \cite{ruthphd}, we showed:

\begin{thm} \label{vnreppsl2r} (Theorem 35 in \cite{ruthphd}.)
Let $G=PSL(2,\R)$, and let $\Gamma$ and $\Lambda$ be lattices in $G$.  Then there exists a representation of the II$_1$ factor $R\Gamma$ on a subspace of $L^2(\Lambda \backslash G)$.  
\end{thm}

To prove this result, we started with the representation of $R\Gamma$ on a discrete series representation in $L^2(G)$ given by:

\begin{thm} \label{vndimthm}
(Theorem 3.3.2 in \cite{ghj}.) 
Let $G$ be a connected semi-simple real Lie group with Haar measure
$dg$, let $\Gamma$ be a discrete subgroup in $G$, and let $\pi : G \rightarrow \UU(\HH)$ be
an irreducible representation in the discrete series. Assume that every non-trivial conjugacy class of $\Gamma$ has infinitely many elements, so that the group von Neumann algebra $R\Gamma$ is a II$_1$ factor. Then the representation $\pi$ of $G$ restricted to $\Gamma$ extends to a representation of the II$_1$ factor $R\Gamma$, with von Neumann dimension
\begin{align} \label{vndimform}
\emph{dim}_{R\Gamma}(\HH)=\emph{covol}(\Gamma)\cdot \emph{d}_\pi
\end{align}
where $\emph{d}_\pi$ is the formal dimension of $\pi$.
\end{thm} 

For background on Theorem \ref{vndimthm} and each piece of formula (\ref{vndimform}), including exposition on discrete series representations and their formal dimensions, on lattices in algebraic groups and their covolumes, and on II$_1$ factors, see \cite{ruthprod}, where we computed von Neumann dimensions in the setting of $PGL(2,F)$ for $F$ a non-archimedean local field,  extending Theorem 61 in \cite{ruthphd}.  

\begin{rem} \label{vndimgen}
The proof of Theorem \ref{vndimthm} applies to a wider class of groups than the theorem states:  The only necessary ingredients are a locally compact unimodular group $G$, a lattice $\Gamma < G$ whose non-trivial conjugacy classes have infinitely many elements, and an irreducible unitary representation of $G$ that is a subrepresentation of the right regular representation of $G$ on $L^2(G)$.  

For example, $GL(n,F)$, $n \geq 2$, $F$ a non-archimedean local field of characteristic $0$ and residue field characteristic relatively prime to $n$, has irreducible unitary representations with matrix coefficients that are compactly-supported modulo the center of $GL(n,F)$, hence square-integrable modulo the center of $GL(n,F)$ (see the introduction of \cite{cms} for more about these representations); and if we consider only those with trivial central character, we obtain irreducible unitary representations of $PGL(n,F)$ with square-integrable matrix coefficients, which are (by Theorem 16.2 in \cite{rob}) equivalent to subrepresentations of the right regular representation of $PGL(n,F)$ on $L^2(PGL(n,F))$.  

For another example, the group $PSL(2,\R) \times \cdots \times PSL(2,\R)$ has representations equivalent to subprepresentations of the right regular representation of \par \noindent
$PSL(2,\R) \times \cdots \times PSL(2,\R)$ on $L^2( PSL(2,\R) \times \cdots \times PSL(2,\R) )$: tensor products of discrete series representations of $PSL(2,\R)$ (by which we mean discrete series representation of $SL(2,\R)$ with trivial central character).  
\end{rem}

The left side of the diagram below illustrates the situation in Theorem \ref{vndimthm}, while the right side of the diagram illustrates the situation in Theorem \ref{vnreppsl2r}, and we use the notation $\HH$ for a discrete series representation, $\alpha$ for a basis vector in $L^2(G / \Gamma)$,  $\beta$ for a basis vector in $l^2(\Lambda)$, $\Phi$ for the intwertwiner between $G$-representations, and $\bar{\Phi}$ for the extension to $R\Gamma$ of the restriction of $\Phi$ to $\Gamma$. 

\[
\begin{tikzcd}
l^2(\Gamma)  \arrow[d, equal] \\
\C \alpha \otimes l^2(\Gamma)  \arrow[d, hookrightarrow] \\
L^2(G / \Gamma) \otimes l^2(\Gamma) \arrow[r, equal] & L^2(G)  \arrow[r, equal] & l^2(\Lambda) \otimes L^2(\Lambda \backslash G) \\
& \HH  \arrow[u, hookrightarrow] \arrow[ddr, rightarrow, "\bar{\Phi}"', bend right=20] & l^2(\Lambda) \otimes  \HH \arrow[u, hookrightarrow] \\
& & \C \beta \otimes  \HH \arrow[u, hookrightarrow] \\
& & \HH \arrow[u, equal] \\
\end{tikzcd}
\]

To obtain the representation of $R\Gamma$ on a discrete series representation of $G$ in $L^2(\Lambda \backslash G)$, instead of in $L^2(G)$, we had to take care of the right-hand side of the diagram above --- that is, we needed to show that some discrete series representation of $G$ occurs in $L^2(\Lambda \backslash G)$.  Since we were working in the setting $PSL(2,\R)$, we used the fact that the multiplicity of a discrete series representation of weight $k \geq 2$ is equal to the dimension of the space of cusp forms of weight $k$ (see \cite{gelb} Theorem 2.10; the proof given there works for all lattices in $SL(2,\R)$, not just $SL(2,\Z)$).  Then formulas for dimensions of spaces of cusp forms (see \cite{miy} Theorem 2.5.2) show that for $k$ high enough, the dimension of the space of cusp forms of weight $k$ for $\Lambda$ is non-zero. (We only considered $k=2, 4, 6, \ldots$ because the even $k$ correspond to the discrete series representations of $SL(2,\R)$ that factor through $PSL(2,\R)$.)  

In the proof of Observation \ref{obsv}, the picture is the same, and a non-zero intertwiner $\Phi$ is guaranteed by the pointwise limit multiplicity property, which we now describe.  

For background on the general limit multiplicity property, see \cite{mulleraastf}, Section 4.  In this note, we are concerned only with the ``pointwise'' limit multiplicity property, since it is all we need to guarantee occurrence of discrete series representations.  If a sequence of lattices has the limit multiplicity property, then it has the pointwise limit multiplicity property. 

Let $\Pi(G)$ denote the unitary dual of $G$ (set of irreducible unitary representations of $G$), 
and let $\Pi(G)_{\text d}$ denote the discrete part of the unitary dual (the representations equivalent to subrepresentations of the right regular representation on $L^2(G)$, \textit{i.e}.\ discrete series representations).

Each $\pi \in \Pi(G)_{\text d}$ has a \textit{formal dimension} $\text d _\pi$ defined by 

\begin{align*}
    \int_G \langle \pi(g) u, v \rangle \overline{\langle \pi(g) u', v'\rangle} dg = \text{d}_\pi^{-1} 
    \langle u, u'\rangle \overline{\langle v,v' \rangle}
\end{align*}

\begin{defn}
We say that a sequence of lattices $( \Lambda_n )$ in a locally compact group $G$ has the \textit{pointwise limit multiplicity property} if 
\begin{align*}
    \lim_{n \rightarrow \infty} \frac{\text{mult} (\pi, \Lambda_n )}{\text{vol}(\Lambda_n \backslash G)} =
    \begin{cases}
        \text d _\pi & \pi \in \Pi(G)_{\text d} \\
        0 & \pi \in \Pi(G) \backslash \Pi(G)_{\text d}
    \end{cases}
\end{align*}
where $\text{mult} (\pi, \Lambda_n )$ denotes the multiplicity of $\pi$ in $L^2(\Lambda_n \backslash G)$.
\end{defn}

(Both $\text d _\pi$ and $\text{vol}(\Lambda_n \backslash G)$ depend on the choice of Haar measure on $G$.)

So, if $\pi$ is a discrete series representation of $G$, and if $( \Lambda_n )$ is a sequence of lattices in $G$ having the pointwise limit multiplicity property, then $\pi$ will ``eventually'' occur in $L^2(\Lambda_n \backslash G)$.

\section{Proof of Observation \ref{obsv}} \label{proof}

\begin{proof}
The generalization of Theorem \ref{vndimthm} discussed in Remark \ref{vndimgen} guarantees that a discrete series representation of $G$ in $L^2(G)$ restricted to $\Gamma$ extends to a representation of $R\Gamma$.  (Now $R\Gamma$ is represented on a subspace $L^2(G)$.)

The sequence $(\Lambda_n)$ in Case 1 has the pointwise limit multiplicity property, by Corollary 1.5 and Remark 1.6 of \cite{gellev}.  

The sequence $(\Lambda_n)$ in Case 2 has the pointwise limit multiplicity property, by Theorem 1 of \cite{matz}.  

The pointwise limit multiplicity property guarantees that the discrete series representation of $G$ occurs in $L^2(\Lambda_i \backslash G)$ for some $\Lambda_i$ in $(\Lambda_n)$.  Let $\Phi$ denote the intertwiner between the representation of $G$ on a subspace of $L^2(G)$ and on a subspace of $L^2(\Lambda_i \backslash G)$.  Then $\Phi$ restricts to an intertwiner between representations of $\Gamma$, and this extends to an intertwiner $\bar{\Phi}$ between representations of $R\Gamma$; see Lemmma 30 of \cite{ruthphd} for the complete argument.  (Now $R\Gamma$ is represented on a subspace $L^2(\Lambda_i \backslash G)$.)
\end{proof}

\begin{rem}
In \cite{matz}, the limit multiplicity property was proven to hold using Sauvageot's density principle and the trace formula.  In \cite{gellev}, the limit multiplicity property was proven to hold using Sauvageot's density principle and convergence of invariant random subgroups, similarly to a theorem in \cite{7s} that connects Benjamini--Schramm convergence (a geometric formulation of convergence of invariant random subgroups) to convergence of relative Plancherel measures. We note that \cite{lev} proved Benjamini--Schramm convergence for sequences of lattices including the sequences considered in \cite{matz}, but the theorems on convergence of relative Plancherel measures in \cite{7s} and \cite{gellev} cannot be applied to give the limit multiplicity property, as the sequences in \cite{matz} consist of non-uniform lattices, while the sequences in \cite{7s} and \cite{gellev} consist of uniform lattices.  
\end{rem}

\begin{que}
All the representations of II$_1$ factors $R\Gamma$ in this note start from restricting discrete series representations of $G$ to $\Gamma$.  Is it possible to represent $R\Gamma$ on a principal series representation or a complementary series representation of $G$, rather than on a discrete series representation of $G$?  Such a representation would not come from restricting a representation of $G$ to $\Gamma$.  
\end{que}

\section*{Acknowledgments} We thank Alain Valette for many enlightening correspondences.  

\bibliographystyle{amsalpha}
\bibliography{ref}

\newcommand{\etalchar}[1]{$^{#1}$}
\providecommand{\bysame}{\leavevmode\hbox to3em{\hrulefill}\thinspace}
\providecommand{\MR}{\relax\ifhmode\unskip\space\fi MR }
\providecommand{\MRhref}[2]{%
  \href{http://www.ams.org/mathscinet-getitem?mr=#1}{#2}
}
\providecommand{\href}[2]{#2}
\begin{thebibliography}{ABB{\etalchar{+}}17}

\bibitem[ABB{\etalchar{+}}17]{7s}
Miklos Abert, Nicolas Bergeron, Ian Biringer, Tsachik Gelander, Nikolay
  Nikolov, Jean Raimbault, and Iddo Samet, \emph{On the growth of
  {$L^2$}-invariants for sequences of lattices in {L}ie groups}, Ann. of Math.
  (2) \textbf{185} (2017), no.~3, 711--790. \MR{3664810}

\bibitem[AS77]{atischmid}
Michael Atiyah and Wilfried Schmid, \emph{A geometric construction of the
  discrete series for semisimple {L}ie groups}, Invent. Math. \textbf{42}
  (1977), 1--62. \MR{0463358}

\bibitem[Ati76]{ati}
M.~F. Atiyah, \emph{Elliptic operators, discrete groups and von {N}eumann
  algebras}, 43--72. Ast\'{e}risque, No. 32--33. \MR{0420729}

\bibitem[CMS90]{cms}
Lawrence Corwin, Allen Moy, and Paul~J. Sally, Jr., \emph{Degrees and formal
  degrees for division algebras and {${\rm GL}_n$} over a {$p$}-adic field},
  Pacific J. Math. \textbf{141} (1990), no.~1, 21--45. \MR{1028263}

\bibitem[GdlHJ89]{ghj}
Frederick~M. Goodman, Pierre de~la Harpe, and Vaughan F.~R. Jones,
  \emph{Coxeter graphs and towers of algebras}, Mathematical Sciences Research
  Institute Publications, vol.~14, Springer-Verlag, New York, 1989. \MR{999799}

\bibitem[Gel75]{gelb}
Stephen~S. Gelbart, \emph{Automorphic forms on ad\`ele groups}, Princeton
  University Press, Princeton, N.J.; University of Tokyo Press, Tokyo, 1975,
  Annals of Mathematics Studies, No. 83. \MR{0379375}

\bibitem[GL18]{gellev}
Tsachik Gelander and Arie Levit, \emph{Invariant random subgroups over
  non-{A}rchimedean local fields}, Math. Ann. \textbf{372} (2018), no.~3-4,
  1503--1544. \MR{3880306}

\bibitem[Lev17]{lev}
Arie Levit, \emph{On benjamini--schramm limits of congruence subgroups},
  \url{arXiv:1705.04200}, 2017.

\bibitem[M\"16]{mulleraastf}
Werner M\"{u}ller, \emph{Asymptotics of automorphic spectra and the trace
  formula}, Families of automorphic forms and the trace formula, Simons Symp.,
  Springer, [Cham], 2016, pp.~477--529. \MR{3675174}

\bibitem[Mat19]{matz}
Jasmin Matz, \emph{Limit multiplicities for {$SL_2(\OO_F)$} in {$SL_2(\R^{r_1}
  \oplus \C^{r_2})$}}, \url{doi: 10.4171/GGD/507} (to appear in print), 2019,
  \textit{Groups Geom. Dyn.} Electronically published on May 7, 2019.

\bibitem[Miy06]{miy}
Toshitsune Miyake, \emph{Modular forms}, english ed., Springer Monographs in
  Mathematics, Springer-Verlag, Berlin, 2006, Translated from the 1976 Japanese
  original by Yoshitaka Maeda. \MR{2194815}

\bibitem[MVN36]{roo1}
F.~J. Murray and J.~Von~Neumann, \emph{On rings of operators}, Ann. of Math.
  (2) \textbf{37} (1936), no.~1, 116--229. \MR{1503275}

\bibitem[Rob83]{rob}
Alain Robert, \emph{Introduction to the representation theory of compact and
  locally compact groups}, London Mathematical Society Lecture Note Series,
  vol.~80, Cambridge University Press, Cambridge-New York, 1983. \MR{690955}

\bibitem[Rut18]{ruthphd}
Lauren~C. Ruth, \emph{Two new settings for examples of von {N}eumann
  dimension}, \url{arxiv:1811.11749}, 2018, dissertation.

\bibitem[Rut19]{ruthprod}
\bysame, \emph{The product of lattice covolume and discrete series formal
  dimension: {$\fp$}-adic {$GL(2)$}}, \url{arXiv:1901.11501}, 2019.

\end{thebibliography}

\end{document}